\documentclass[10pt,a4paper]{amsart}

 \usepackage[latin1]{inputenc}
 \usepackage[T1]{fontenc}
\usepackage{t1enc}
\usepackage{xy}
 \usepackage{eepic,epic}
 \usepackage{mathrsfs}
 \usepackage{graphics,psfrag}
\usepackage{graphicx} 
\xyoption{all}

\def\Ddots{\mathinner{\mkern1mu\raise1pt\hbox{.}\mkern2mu\raise4pt\hbox{.}\mkern
 2mu
   \raise7pt\vbox{\kern7pt{\hbox{.}}}\mkern1mu}}

\theoremstyle{plain}
\newtheorem{theo}{Theorem}
\newtheorem{lemm}{Lemma}
\newtheorem{prop}{Proposition}
\newtheorem{coro}{Corollary}

\theoremstyle{definition}
\newtheorem{exem}{Example}
\newtheorem{defi}{Definition}

\theoremstyle{remark}
\newtheorem{rema}{Remark}

\title{A chevalley formula in equivariant $K$-theory}

\author[Matthieu Willems ]{Matthieu Willems}

\address{Universit{\'e} de Cergy-Pontoise \\
D{\'e}partement de math{\'e}matiques \\ 
Site de Saint-Martin \\
2 rue Adolphe Chauvin \\
  95302 Cergy-Pontoise Cedex \\
  France}

\email{matthieu.willems@polytechnique.org}

\begin{document}

\subjclass{$19$L$47$, $14$M$15$ }


\maketitle

\tableofcontents

\begin{abstract}

The aim of this paper is to give a recursive formula to compute the
product of a line bundle with the structure sheaf of a Schubert
variety in the equivariant $K$-theory of a flag variety.

\end{abstract}

\section{Introduction}

Let $G$ be a complex semi-simple connected group, $B
\subset G$ a Borel subgroup of $G$, and $H \subset B$ a
maximal torus of $B$. We denote by $R[H]$ the ring of 
representations of $H$ and $X=G/B$ the flag variety of $G$.
The $H$-equivariant $K$-theory $K(H,X)$ of $X$ has a $R[H]$-basis 
$[\mathcal{O}_{\overline{X}_w}]^H$ indexed by $W=N_G(H)/H$ the Weyl
group of $G$, where $[\mathcal{O}_{\overline{X}_w}]^H$ is the class of
the structure sheaf of the Schubert variety $\overline{X}_w$. The
Schubert variety $\overline{X}_w \subset X$ is the closure of the
$B$-orbit of $w \in W$. 
Let $\mathfrak{h}$ be the Lie algebra of $H$, we denote by 
$\mathfrak{h}_{\mathbb{Z}}^* \subset \mathfrak{h}^*$ the weight
lattice which is identified canonically with the group of characters
of $H$. Then for all $\lambda \in \mathfrak{h}_{\mathbb{Z}}^*$, we 
denote by $\mathcal{L}_{\lambda}^X$ the canonical line bundle over
$X$.   
The torus $H$ acts on $\mathcal{L}_{\lambda}^X$ and it defines a class
$[\mathcal{L}_{\lambda}^X]^H$ in $K(H,X)$. In fact,
$K(H,X)$ is generated by these line bundles and then we get a
presentation by generators and relations of $K(H,X)$.
If we want to understand the link between the basis
$\{ [\mathcal{O}_{\overline{X}_w}]^H \}_{w \in W}$ and this
presentation, it is interesting to find a "Giambelli formula" which
expresses $[\mathcal{O}_{\overline{X}_w}]^H$ in terms of
$\{[\mathcal{L}_{\lambda}^X]^H \}_{\lambda \in
  \mathfrak{h}_{\mathbb{Z}}^*}$ and a "Chevalley formula", i.e. to
find the coefficients $q_{w, v}^{\lambda} \in R[H]$ satisfying : 
$$[\mathcal{L}_{\lambda}^X]^H [\mathcal{O}_{\overline{X}_w}]^H= 
\sum_{v \in W}q_{w, v}^{\lambda} [\mathcal{O}_{\overline{X}_v}]^H.
$$
Such a formula has been known for a long time in cohomology (see
\cite{che} in ordinary cohomology and \cite{kkc} in equivariant
cohomology).
In \cite{pittie} Pittie and Ram give a Chevalley formula in ordinary
$K$-theory for a dominant weight $\lambda$ by using L-S paths. 
Littelmann and Seshadri generalize this formula to $H$-equivariant
$K$-theory in \cite{littel}.  
Such a formula was first given in the case $G=SL(n,\mathbb{C})$ by
Fulton and Lascoux in \cite{ful-las} by using "tableaux" of shape
$\lambda$.  
The aim of this paper is to find a practical algorithm to compute
these coefficients $q_{w, v}^{\lambda}$. Our formula is valid for all
weights (even for non dominant weights).

\medskip

We follow the same method as in \cite{mw3} to find restrictions to
fixed points in equivariant cohomology and K-theory. First we describe
a $R[H]$-basis of the $H$-equivariant $K$-theory of a Bott-Samelson
variety $\Gamma$ and we decompose the class of a line bundle
$[\mathcal{L}_{\lambda}^{\Gamma}]^H$ in this basis. To find this formula,
we use the structure of iterated fibrations with fiber $\mathbb{C}P^1$
of Bott-Samelson varieties. Then we use the standard map $g : \Gamma
\rightarrow X$ to deduce a Chevalley formula in $K(H,X)$. In
\cite{duan3} and \cite{duan4} Haibao Duan also used Bott-Samelson
varieties to find formulas in Schubert calculus and we used this idea
in \cite{mw4} and \cite{mw5} to find similar formulas in the
equivariant setting. In these two papers we study Bott towers i.e. all
varieties which have a structure of iterated fibrations with fiber
$\mathbb{C}P^1$.

 \bigskip

In Section~\ref{notations} we recall basic definitions on semi-simple 
groups and their flag varieties.

In Section~\ref{defBS}, we recall the definition of  the Bott-Samelson
variety associated to a sequence of simple roots and we define a cell
decomposition of this variety. For more details on this section, see
\cite{livrekumar} and \cite{mw4}. 

In Section~\ref{ktheory}, we recall the definition of the
$H$-equivariant $K$-theory of an algebraic $H$-variety and we
introduce the notion of restriction to fixed points which will be the
main tool of our proofs.

\medskip

In Section~\ref{sectionKtheoryBott}, we construct a $R[H]$-basis of the 
$H$-equivariant $K$-theory of a Bott-Samelson variety $\Gamma$ and 
for all $\lambda \in \mathfrak{h}_{\mathbb{Z}}^*$, we decompose the 
line bundle $\mathcal{L}_{\lambda}^{\Gamma}$ in this basis 
(Theorem \ref{linebundleBS}).

\medskip

In Section~\ref{sectionvd}, if $g : \Gamma \rightarrow X$ is the standard 
map from a Bott-Samelson variety $\Gamma$ to the flag variety $X$, we
describe the morphism  $g_*$ induced in $K$-theory (Theorem \ref{g*}) and we 
deduce from this result the main theorem of this paper (Theorem \ref{main}) 
which gives a Chevalley formula in equivariant $K$-theory.

\medskip

In Section~\ref{sectionordinary}, we restrict our calculations to ordinary 
$K$-theory (Theorem \ref{mainordinary}).

\bigskip

I am very grateful to Michel Brion for suggesting me using my
results to find a Chevalley formula in equivariant $K$-theory.

\section{Preliminaries and notation} \label{notations}

\subsection{Root system}

Let $G$ be a  connected and simply connected complex semi-simple
group of rank $r$. 
We denote by $e$ the neutral element of $G$. 
Let $B \subset G$ be a Borel subgroup of $G$ and $H \subset B$ the
Cartan subgroup of $B$.
We denote by $\mathfrak{h} \subset \mathfrak{b} \subset \mathfrak{g}$
the Lie algebras of $H$, $B$ and $G$. 
\medskip

We choose a system of simple roots $\mathfrak{\pi} = \{\alpha_{i}\}_{1
 \leq i \leq  r} \subset \mathfrak{h}^*$ and simple coroots 
$\mathfrak{\pi^{\vee}} = \{\alpha_i^{\vee}\}_{1 \leq i \leq r}  \subset
 \mathfrak{h}$, such that  
$$\mathfrak{b}=\mathfrak{h} \oplus \bigoplus_{\alpha \in
  \Delta_{+}}\mathfrak{g}_{\alpha} \,\, , \,\,\,\,\, {\rm and \,\,\,\,} 
\mathfrak{g}=\mathfrak{h} \oplus \bigoplus_{\alpha \in
  \Delta_{+}}(\mathfrak{g}_{\alpha} \oplus \mathfrak{g}_{-\alpha}),$$
where for $\lambda \in \mathfrak{h}^*$, $\mathfrak{g}_{\lambda} = \{ x
  \in \mathfrak{g} \: {\rm \, such \: that \, }\: [h, x]=\lambda(h)x,
  \forall h \in \mathfrak{h} \}$, and where we define $\Delta_{+}$ by
   $\Delta_{+} = \{ \alpha \in \sum_{i=1}^{r}\mathbb{N}\alpha_{i} \:
   {\rm \, such \: that \,}  \: \alpha \neq 0
 \:{\rm and } \: \mathfrak{g}_{\alpha} \neq 0 \}$. We set
   $\Delta=\Delta_{+} \cup \Delta_{-}$ where $\Delta_{-} =
   -\Delta_{+}$. We call $\Delta_{+}$ 
 (respectively $\Delta_{-}$) the set of positive roots (respectively
 negative).

\medskip

We associate to $(\mathfrak{g}, \mathfrak{h})$ the Weyl group
$W\subset { \rm Aut}(\mathfrak{h}^*)$ generated by the simple
reflections  $\{s_{i}\}_{1  \leq i \leq r}$ defined by 
$$ \forall \lambda \in \mathfrak{h}^*, \,\,
 s_{i}(\lambda)=\lambda-\lambda(\alpha^{\vee}_{i})\alpha_{i} .$$
By dualizing, we get an action of $W$ on $\mathfrak{h}$.

If we denote by $S$ the set of simple reflections, the couple $(W,S)$
is a Coxeter system. Thus we have a notion of Bruhat order denoted by  
$u \leq v$ and a notion of length
denoted by $l(w) \in \mathbb{N}$. We  denote by $1$ the neutral
element of $W$. 

We have $\Delta=W\pi$, and for $\beta = w\alpha_{i}
\in \Delta^{+}$, we set $s_{\beta}=ws_{i}w^{-1} \in W$ (which does not  
depend on the couple $(w,\alpha_{i})$ satisfying $\beta =
w\alpha_{i}$), and $\beta^{\vee} = w\alpha_i^{\vee} \in \mathfrak{h}$.

\medskip

We define the fundamental weights $\rho_{i} \in \mathfrak{h}^*$ 
($1 \leq i \leq r$) by 
$$\rho_{i}(\alpha^{\vee}_j)=\delta_{i, j}, { \rm \, for \,\, all \,\, } 
1 \leq i,j \leq r,$$ 
and the weight lattice $\mathfrak{h}_{\mathbb{Z}}^*$ by 
$$\mathfrak{h}_{\mathbb{Z}}^* = \oplus_{1 \leq i \leq
  r}\mathbb{Z} \rho_i \subset  \mathfrak{h}^*.  
$$

\subsection{Flag varieties}

Let $N_{G}(H)$ be the normalizer of $H$ in $G$, the quotient group
$N_{G}(H)/H$ can be identified to $W$. We set $X=G/B$. It is a flag
variety. 
 The group $G$ acts on $X$ by multiplication on the left. This action
 yields  an action of $B$ and $H$ on $X$. The set of fixed points of
 this action of $H$ on $X$ can be  identified to $W$.
 For $w \in W$, we define $C(w)=B \cup BwB$ and for all simple roots
 $\alpha$, we define the subgroup $P_{\alpha}$ of $G$ by
 $P_{\alpha}=C(s_{\alpha})$. We have the Bruhat decomposition of
 $G=\bigsqcup_{w \in W}BwB$ and if we define $X_{w}=BwB/B$, then
 $X=\bigsqcup_{w \in W}X_{w}$. For all $w \in W$, the Schubert cell
 $X_{w}$ is isomorphic to $\mathbb{C}^{l(w)}$. Thus we get an
 $H$-equivariant cell decomposition of $X$ where all cells have even 
 real dimension.  

\medskip

For all $w \in W$, the Schubert variety $\overline{X}_{w}$ is the
closure of the cell $X_{w}$. 
It is an irreducible $H$-equivariant subvariety of $X$ of complex
dimension $l(w)$. In general Schubert varieties are not smooth. 
For all $w\in W$, we have the decomposition 
$$\overline{X}_{w}=\bigsqcup_{w' \leq w }X_{w'}.$$

\subsection{The monoid \underline{$W$} } \label{monoid}

We define the monoid $\underline{W}$ as the monoid generated by the
 elements $\{\underline{s}_{i}\}_{1 \leq i \leq r}$ with the relations
$\underline{s}_{i}^2=\underline{s}_{i}$ and the braid relations of $W$ :  
$$ \left\{ \begin{array}{cc}
     \underline{s}_{i}^2=\underline{s}_{i} & \\ 
  \underbrace{\underline{s}_{i}\, \underline{s}_{j} \cdots
     }_{m_{i,j}\, { \rm terms}}=   \underbrace{\underline{s}_{j}\,  
\underline{s}_{i} \cdots }_{m_{i,j} \, { \rm      terms}  }
  & { \rm if }\,  m_{i,j}<\infty \, ,
   \end{array} \right.$$
where $m_{i,j}$ is the order of  $s_is_j$ in $W$.

\medskip

We denote by $T : W \rightarrow \underline{W}$ the bijection defined
by $T(w)=\underline{s}_{i_{1}} \cdots 
\underline{s}_{i_{l}}$ if $w=s_{i_{1}} \cdots s_{i_{l}}$ is a reduced  
decomposition of $w$ (i.e. $l=l(w)$).




\section{Bott-Samelson varieties} \label{defBS}

Let $N \geq 1$ be a positive integer. We use the notation of 
Section~\ref{notations}.

\subsection{Definition} \label{def}

Let $\mu_1$, \ldots , $\mu_N$ be a sequence of $N$ simple roots
(repetitions may occur). We define 
$$\Gamma(\mu_{1}, \ldots ,\mu_{N})=P_{\mu_1} \times_{B} 
P_{\mu_2} \times_{B} \cdots \times_{B} P_{\mu_N}/B,$$
as the space of orbits of $B^N$ acting on $P_{\mu_1} \times P_{\mu_2} 
\times \cdots \times P_{\mu_N}$ by
$$(g_{1}, g_{2}, \ldots , g_{N}) (b_{1}, b_{2}, \ldots , b_{N}) = 
(g_{1}b_{1},b_{1}^{-1} g_{2}b_{2}, \ldots ,b_{N-1}^{-1} g_{N}b_{N}),\,
 b_{i} \in B, \, g_{i} \in P_{\mu_i}.$$

\medskip

It is an irreducible complex projective variety of dimension $N$. 
We denote by $[g_{1}, g_{2}, \ldots , g_{N}]$ the class of $(g_{1},
g_{2}, \ldots , g_{N})$ in $\Gamma(\mu_{1}, \ldots ,\mu_{N})$ and by  
$g_{\mu_{i}} \in P_{\mu_i}$ a representative of the reflection 
of $N_{P_{\mu_i}}(H)/H \simeq \mathbb{Z}/2\mathbb{Z}$.

\bigskip

We define a left action of $B$ on $\Gamma(\mu_1, \ldots, \mu_N)$ by  
$$b[g_{1}, g_{2}, \ldots ,g_{N}]=[bg_{1}, g_{2}, \ldots
,g_{N}],\hspace{0,1 cm} b \in B, \hspace{0,1 cm}  g_{i} \in
P_{\mu_i}.$$ 

By restricting this action to $H$, we get an action of $H$ on
$\Gamma(\mu_1, \ldots, \mu_N)$. 

In the following two sections we denote $\Gamma(\mu_1, \ldots,
\mu_N)$ by  $\Gamma$.

\subsection{Cell decomposition}\label{decomposition}

 For $\epsilon \in \{0,1\}^N$, we denote by $\Gamma_{\epsilon}
\subset \Gamma$ the set of classes $[g_{1},
g_{2}, \ldots , g_{N}]$ satisfying for all integers  $1 \leq i \leq N$ 
$$ g_{i} \in B
  {\,\,\, \rm if} \hspace{0,15 cm} \epsilon_{i} =0 \, ,
  \,\,\,\,\,\,\,\,  g_{i} \notin B
  {\,\,\, \rm if }\hspace{0,15 cm} \epsilon_{i} =1.
$$

\medskip

For $\epsilon =(\epsilon_1, \epsilon_2, \ldots, \epsilon_N) \in
\{0,1\}^N$, we denote by $l(\epsilon)$ the cardinal of
$\{ 1 \leq i \leq N, \epsilon_i=1\}$. It is called the length of
$\epsilon$. We define a partial order on $\{0,1\}^N$ by 
 $$\epsilon \leq \epsilon' \Leftrightarrow (\forall 1 \leq i \leq N,  
\epsilon_i = 1 \Rightarrow  \epsilon'_i=1).$$

\bigskip

The following proposition is obvious.

\begin{prop} \label{structureBS}
    
\indent

\begin{enumerate}

\item[$(i)$] For all $\epsilon \in \{0,1\}^N$, $\Gamma_{\epsilon}$ 
 is a complex affine space of dimension $l(\epsilon)$ 
 which is invariant under the action of $B$, 
and this action induces a linear action of the torus $H$ on
 $\Gamma_{\epsilon}$.  

\item[$(ii)$] For all $\epsilon \in \{0,1\}^N$,
  $\overline{\Gamma}_{\epsilon} =  
\coprod_{\epsilon' \leq \epsilon} \Gamma_{\epsilon'}$.

\item[$(iii)$] $\Gamma = \coprod_{\epsilon \in
  \{0,1\}^N} \Gamma_{\epsilon}$.

\item[$(iv)$] For all $\epsilon  \in \{0,1\}^N$,
  $\overline{\Gamma}_{\epsilon}$ can be identified
with the Bott-Samelson variety $\Gamma(\mu_{i}, \, 
 \epsilon_i=1)$ and is 
an irreducible smooth subvariety of $\Gamma$.

\end{enumerate}   

\end{prop}

 For $\epsilon \in \{0,1\}^N$ and $1 \leq i\leq N$, we define 
$$\displaystyle{v_{i}(\epsilon) =\prod_{ {\tiny \begin{array}{cc}
 1\leq j \leq i, \\   \epsilon_j=1 \end{array}} 
  }\!\!\!\!\!\!s_{\mu_{j}}}  \in W,$$
where, by convention, $\prod_{\emptyset}=1$.  We set
 $v(\epsilon)=v_{N}(\epsilon) \in W$.

\medskip

Moreover we define the root  
$\alpha_{i}(\epsilon) \in \Delta$ by 
$$\alpha_{i}(\epsilon)=v_{i}(\epsilon)\mu_{i}.$$  

\bigskip

Let $\Gamma^H$ be the set of fixed points of the action of $H$ on
$\Gamma$, we can identify $\Gamma^H$ with $\{0,1\}^N$ thanks 
to the following lemma.

\begin{lemm} \label{fixedpointsBS}

\indent 
    
\begin{enumerate}

\item[$(i)$]  $\Gamma^H \simeq \displaystyle \prod_{1 \leq i \leq N}
  N_{P_{\mu_i}}(H)/H \simeq \prod_{1
  \leq i \leq N}  \{ e, g_{\mu_{i}} \} \simeq \{0,1\}^N,$
 where we identify $e$ with $0$ and $g_{\mu_{i}}$ with $1$.

\item[$(ii)$] For all $\epsilon \in \{0,1\}^N$, $\Gamma_{\epsilon}$
  is the $B$-orbit of $\epsilon \in \Gamma^H$.

\item[$(iii)$] For $(\epsilon,\epsilon') \in (\{0,1\}^N)^2$ 
$$ \epsilon \in \overline{\Gamma}_{\epsilon'} \Leftrightarrow \epsilon 
\leq \epsilon',$$ 
and if we denote by $T_{\epsilon'}^{\epsilon}$ the tangent space to 
$\overline{\Gamma}_{\epsilon'}$ at $\epsilon$, then the weights of the
representation of  $H$ in $T_{\epsilon'}^{\epsilon}$ are
$\{ -\alpha_i(\epsilon) \}_{i, \, \epsilon'_i=1 }$.

\end{enumerate}

\end{lemm}

\subsection{Fibrations of Bott-Samelson varieties} \label{fibrations}

For all $2\leq k \leq N$, let denote by $\pi_k :  \Gamma(\mu_1,
\ldots,
\mu_k) \rightarrow  \Gamma(\mu_1,  \ldots, \mu_{k-1})$ the projection
defined by 
$$ \pi_k([g_1, \ldots, g_k])=[g_1, \ldots, g_{k-1}].
$$ 
If we denote by $\pi_1 : \Gamma(\mu_1) \rightarrow \{ \rm point \}$
the trivial projection, we get the following diagram : 

 $$ \begin{array}{c}
   \Gamma(\mu_1,  \ldots, \mu_N) \\ 
  \downarrow \pi_N   \\
 \Gamma(\mu_1,  \ldots, \mu_{N-1}) \\
 \,\,\,\, \,\,\, \downarrow \pi_{N-1}   \\
\!\!\!\!\!\!\!\! \vdots  \\
 \downarrow \pi_{3}   \\
  \Gamma(\mu_1,\mu_2) \\
   \downarrow \pi_2 \\
   \Gamma(\mu_1) \simeq \mathbb{C}P^1   \\
  \,\, \downarrow \pi_1   \\
 \{ \rm point\}    
\end{array}$$

where each projection $\pi_k$ is a fibration with fiber
$\mathbb{C}P^1$ (see \cite{mw4} for more details).

\subsection{Line bundles} \label{sectionlinebundlesBS}

We denote by $X(H)$ the group of characters of $H$.
For all integral weights $\lambda \in \mathfrak{h}_{\mathbb{Z}}^*$, 
we denote by $e^{\lambda} : H \rightarrow S^1$ the
corresponding character. This way we get an isomorphism between the 
additive group $\mathfrak{h}_{\mathbb{Z}}^*$ and $X(H)$.

\medskip

Since $H \simeq B/U$, where $U$ is the unipotent radical of $B$, we can 
extend to $B$ all characters $e^{\lambda} \in X(H)$ (in fact 
$X(H) \simeq X(B)$). Then for all $\lambda \in \mathfrak{h}_{\mathbb{Z}}^*$, 
we denote by $\mathcal{L}_{\lambda}^{\Gamma}$ the $B$-equivariant line
bundle over $\Gamma$ defined as the space of orbits of $B^N$ acting on   
$P_{\mu_1} \times P_{\mu_2} \times \cdots \times P_{\mu_N} 
\times \mathbb{C}$ by
$$(g_{1}, g_{2}, \ldots , g_{N}, v) (b_{1}, b_{2}, \ldots , b_{N}) = 
(g_{1}b_{1},b_{1}^{-1} g_{2}b_{2}, \ldots ,b_{N-1}^{-1} g_{N}b_{N}, 
e^{\lambda}(b_N^{-1})v),$$ 
$$ b_{i} \in B, \, g_{i} \in P_{\mu_i}, \, v \in \mathbb{C}.$$

\section{Equivariant $K$-theory} \label{ktheory}

Let $Z$ be a complex algebraic $H$-variety, we denote by $Z^H \subset
Z$ the set of   fixed points of the action of $H$ on $Z$.

 We denote by  $K^0(H,Z)$ the Grothendieck group of $H$-equivariant 
complex  vector bundles of finite rank over $Z$. The tensor
product  of vector bundles defines a product on $K^0(H,Z)$. Since
$K^0(H,{\rm point }) \simeq R[H]$, where  $R[H]=\mathbb{Z}[X(H)] $  
is the representation ring of the torus $H$, we get a  $R[H]$-algebra 
structure on $K^0(H,Z)$.

\medskip

For all $H$-equivariant algebraic maps $g : Z_1 \rightarrow Z_2$ 
we denote by  $g^* : K^0(H,Z_2) \rightarrow K^0(H,Z_1)$
the morphism of $R[H]$-algebras defined by pulling back vector
bundles. In particular the inclusion $Z^H \subset Z$  gives a morphism  
$i_H^*  : K^0(H,Z) \rightarrow K^0(H,Z^H)$ called restriction to fixed
points. If the set of fixed points $Z^H$ is finite, $K^0(H,Z^H) $
can be identified with $F(Z^H; R[H])$ the $R[H]$-algebra of all maps
 $f : Z^H \mapsto R[H]$ and we get a morphism 
 $i_H^* :  K^0(H,Z)  \rightarrow  F(Z^H; R[H])$. Moreover if
 $K^0(H,Z)$ is a free $R[H]$-module, then the morphism $i_H^*$ is
 injective. This is an easy consequence of the localization theorem (see 
 \cite{ginzburg} section $5.10$).

\medskip

If we assume that $Z$ is a complex projective smooth $H$-variety, 
then $K^0(H,Z)$ is isomorphic to $K_{0}(H,Z)$ the Grothendieck group
of $H$-equivariant coherent sheaves on $Z$ (see \cite{ginzburg}
chapter 5). In this case, we identify these two groups and we denote
them by $K(H,Z)$. 

\medskip

For all proper $H$-equivariant morphisms $g : Z_1 \rightarrow Z_2$,  
we denote by  $g_* : K_0(H,Z_1) \rightarrow K_0(H,Z_2)$ the direct
image morphism. 
For all $H$-equivariant subvarieties $Z' \subset Z$, we denote by
$[\mathcal{O}_{Z'}]^H \in K_0(H,Z)$ the class of
$i_*(\mathcal{O}_{Z'})$ where $\mathcal{O}_{Z'}$ is the structure
sheaf of $Z'$ and $i$ is the inclusion of $Z'$ in $Z$.

\section{$K$-theory of Bott-Samelson varieties} \label{sectionKtheoryBott}

We use the notation of Section~\ref{defBS}. Let $N$ be a positive
integer, and let $\mu_1$, $\mu_2$, \ldots, $\mu_N$ be a sequence of
$N$ simple roots. Let $\Gamma$ be the Bott-Samelson variety $\Gamma=
 \Gamma(\mu_1, \mu_2, \ldots , \mu_N)$. For $1 \leq k \leq N$, we
 denote $\Gamma(\mu_1, \ldots , \mu_k)$ by $\Gamma^k$. By convention, 
$\Gamma^0 = \{ \rm point \}$. For $1 \leq k \leq N$, let $\pi_k$ be the
projection $\Gamma^k \rightarrow \Gamma^{k-1}$ defined in
Section~\ref{fibrations}.

We denote by $\displaystyle \Gamma^k=\coprod_{\epsilon \in 
\{0,1\}^k}\Gamma_{\epsilon}^k$ the cell decomposition defined in 
Section~\ref{decomposition}. For all $\epsilon \in \{0,1\}^k$, let 
denote by $\overline{\Gamma}_{\epsilon}^k$ the closure of
$\Gamma_{\epsilon}^k$ in $\Gamma^k$.

\subsection{A basis of the $K$-theory of Bott-Samelson varieties}
\label{basisBS}

For all integers $1 \leq k \leq N$, $\Gamma^k$ is a complex projective
smooth $H$-variety, and then we denote by $K(H,\Gamma^k)$ its
$H$-equivariant $K$-theory.   

For all $\epsilon \in \{0,1\}^k$, we set $\mathcal{O}_{k,\epsilon}^H =
[\mathcal{O}_{\overline{\Gamma}_{\epsilon}^k}]^H \in K(H,\Gamma^k)$,
and for $1 \leq i \leq k$, $\mathcal{O}_{k,i}^H  =\mathcal{O}_{
  k,[^i_k]  }^H$, where for $1 \leq j \leq k$, $[^i_k]_j  =
1-\delta_{i,j}$. 

Since $\Gamma^k = \coprod_{\epsilon \in  \{0,1\}^k} \Gamma_{\epsilon}^k$ 
is a cell decomposition of $\Gamma^k$, the family 
$\{\mathcal{O}_{k,\epsilon}^H \}_{\epsilon \in \{0,1\}^k}$ is a
$R[H]$-basis of the module $K(H,\Gamma^k)$. Moreover, $({\Gamma^k})^H$
is finite and isomorphic to $\{0,1\}^k$. Thus we have the following
proposition.

\begin{prop} \label{propbasektheorie}

\indent

\begin{enumerate}

\item[$(i)$] The $H$-equivariant $K$-theory of $({\Gamma^k})^H$ can be
  identified with $F(\{ 0,1\}^k;R[H])$.

\item[$(ii)$]  $ K(H,\Gamma^k)=\oplus_{\epsilon \in \{0,1\}^k}
R[H] \mathcal{O}_{k,\epsilon}^H$.

\item[$(iii)$] The restriction to fixed points $i_{H}^*$ :
  $K(H,\Gamma^k) \rightarrow F(\{0,1\}^k;R[H])$ is injective.

\end{enumerate}

\end{prop}

For all integers $1 \leq k \leq N$ and all $\lambda \in
\mathfrak{h}_{\mathbb{Z}}^*$, we denote by 
$\mathcal{L}_{\lambda}^k$ the line bundle over $\Gamma^k$ 
defined in Section \ref{sectionlinebundlesBS} and by 
$[\mathcal{L}_{\lambda}^k]^H$ its class in $K(H,\Gamma^k)$.
For $k=0$, $[\mathcal{L}_{\lambda}^0]^H=e^{\lambda} \in R[H]
\simeq K(H,{\rm point})$.

We will decompose $[\mathcal{L}_{\lambda}^k]^H$ in the $R[H]$-basis 
$\{\mathcal{O}_{k,\epsilon}^H \}_{\epsilon \in \{0,1\}^k}$. For this
we need restrictions to fixed points.

\subsection{Restrictions to fixed points}

\begin{prop} \label{restrictions} 
For all integer $1 \leq k \leq N$, $\lambda \in
\mathfrak{h}_{\mathbb{Z}}^*$ and $(\epsilon, \epsilon') \in
(\{0,1\}^k)^2$,

$$i_H^*( [\mathcal{L}_{\lambda}^k]^H)(\epsilon) = 
e^{v(\epsilon)\lambda},$$

$$i_H^*(\mathcal{O}_{k,\epsilon'}^H)(\epsilon)=
\left\{ \begin{array}{ll}
\!\!\!\! \displaystyle{\prod_{1 \leq i \leq k, \epsilon'_i
    =0}(1-e^{-\alpha_{i}(\epsilon)}) }   
 & { \rm if } \hspace{0,15 cm} \epsilon \leq \epsilon', \\ 
  0
 & { \rm  otherwise }.
\end{array}\right.$$

 \end{prop}
 
 \begin{proof}

The first relation is obvious by definition of $\mathcal{L}_{\lambda}^k$.

\medskip

Let us prove the second relation.

If $\epsilon \not\leq \epsilon'$, the fixed point $\epsilon \not\in
\overline{\Gamma}_{\epsilon'}^k$, 
and then by the localization theorem,
$$i_H^*(\mathcal{O}_{k,\epsilon'}^H) (\epsilon)=0.$$ 

If $\epsilon \leq \epsilon'$, $\epsilon \in
\overline{\Gamma}_{\epsilon'}^k$.  
Since $\overline{\Gamma}_{\epsilon'}^k$  and $\Gamma^k$
are smooth, we can use the self-intersection formula (see
\cite{ginzburg} proposition 5.4.10) to get 
$i_H^*(\mathcal{O}_{k,\epsilon'}^H)(\epsilon)=
i_H^*(\lambda(T_{\overline{\Gamma}_{\epsilon'}^k}^*\Gamma^k))(\epsilon)$,
where $T_{\overline{\Gamma}_{\epsilon'}^k}^* \Gamma^k$ is the normal
bundle of $\overline{\Gamma}_{\epsilon'}^k$ in $\Gamma^k$ and for all
vector bundles $V$, $\lambda(V)= \sum_{0 \leq i \leq \dim(V)}(-1)^i
\Lambda^i(V)$.  

\smallskip

Then Lemma \ref{fixedpointsBS} and the relation $\lambda(V_1 \oplus
V_2) =\lambda(V_1)\otimes \lambda(V_2)$ give us
$$i_H^*(\mathcal{O}_{k,\epsilon'}^H)(\epsilon)=
\displaystyle{\prod_{1 \leq i \leq k, \epsilon_i' = 0} \! \! \! \! 
  (1-e^{-\alpha_{i}(\epsilon)}) }.  $$

\end{proof}

Since $i_H^*$ is injective, we deduce the following formula.

\begin{coro} \label{product}
For all integers $1 \leq k \leq N$ and all $\epsilon \in \{0,1\}^k$,
\begin{equation} \mathcal{O}_{k,\epsilon}^H = 
\prod_{1 \leq i \leq k, \epsilon_i = 0} \! \! \!
\mathcal{O}_{k,i}^H.\end{equation} 

\end{coro}

\begin{rema}

This corollary is also a consequence of the fact that
$$\overline{\Gamma}_{\epsilon}^k=\bigcap_{\tiny \begin{array}{c} 1 \leq i
    \leq k, \\ \epsilon_i' = 0 \end{array}} \! \! \! \! \!
\overline{\Gamma}_{[^i_k]}^k$$ is a transversal intersection.

\end{rema}

\subsection{Decomposition of line bundles}

\begin{theo}  \label{induction}

For all integers $1 \leq k \leq N$, and all $\lambda \in
  \mathfrak{h}_{\mathbb{Z}}^*$,

$$[\mathcal{L}_{\lambda}^k]^H=\pi_k^*([\mathcal{L}_{s_{\mu_k} 
\lambda}^{k-1}]^H) + \mathcal{O}_{k,k}^H \pi_k^*([\mathcal{L}_{
\lambda, \mu_k}^{k-1}]^H), 
$$
where $$[\mathcal{L}_{\lambda, \mu_k}^{k-1}]^H = \left\{ \begin{array}{ll}
0  
 & { \rm if } \hspace{0,15 cm}  \lambda(\mu_k^{\vee})=0, \\ 
 
 [\mathcal{L}_{\lambda}^{k-1}]^H 
 +   [\mathcal{L}_{\lambda-\mu_k}^{k-1}]^H+ \cdots +
     [\mathcal{L}_{\lambda-(\lambda(\mu_k^{\vee})-1)\mu_k}^{k-1}]^H
 & { \rm if } \hspace{0,15 cm}  \lambda(\mu_k^{\vee}) >0,  \\
- [\mathcal{L}_{\lambda+\mu_k}^{k-1}]^H- \cdots -
    [ \mathcal{L}_{\lambda-\lambda(\mu_k^{\vee})\mu_k}^{k-1}]^H
 & { \rm if } \hspace{0,15 cm}  \lambda(\mu_k^{\vee}) <0  .
\end{array}\right.
$$

\end{theo}

\begin{proof}

To prove this theorem, we check these relations after restriction to
fixed points.

Let $\epsilon$ be an element of $\{0,1\}^k$. We denote by
$\overline{\epsilon}$ the element of $\{ 0,1\}^{k-1}$ corresponding to
the fixed point $\pi_k(\epsilon)$ in $\Gamma^{k-1}$ (by convention
$\{0,1\}^0 = \{0\}$, and $v(0)=1\in W$). 

If we use Proposition~\ref{restrictions} for $\Gamma^k$ and
$\Gamma^{k-1}$, we find : 
$$i_H^*([\mathcal{L}_{\lambda}^k]^H)(\epsilon)=e^{v(\epsilon)\lambda},$$

$
i_H^*\big(\pi_k^*([\mathcal{L}_{s_{\mu_k} 
\lambda}^{k-1}]^H) + \mathcal{O}_{k,k}^H \pi_k^*([\mathcal{L}_{
\lambda, \mu_k}^{k-1}]^H)\big)(\epsilon)= $

\hspace{5 cm} $ e^{v(\overline{\epsilon})s_{\mu_k}\lambda}+\delta_{\epsilon_k,
  0}(1-e^{-v(\epsilon)\mu_k})  
(i_H^*([\mathcal{L}_{\lambda,
  \mu_k}^{k-1}]^H))(\overline{\epsilon}). 
$

\bigskip

If $\epsilon_k=1$, we have to check $e^{v(\epsilon)\lambda}
=e^{v(\overline{\epsilon})s_{\mu_k}\lambda}$. This is obvious since 
$v(\epsilon)=v(\overline{\epsilon})s_{\mu_k}$.

\medskip

If $\epsilon_k=0$, $v(\epsilon)=v(\overline{\epsilon})$ 
and thus we have to check 
$$ e^{\lambda} = \left\{ \begin{array}{ll}
e^{s_{\mu_k}\lambda}
 & { \rm if } \hspace{0,15 cm}  \lambda(\mu_k^{\vee})=0, \\ 
 e^{s_{\mu_k}\lambda} +(1-e^{-\mu_k})\big[ e^{\lambda}+
 e^{\lambda-\mu_k}+ \cdots +
     e^{\lambda-(\lambda(\mu_k^{\vee})-1)\mu_k}\big]
 & { \rm if } \hspace{0,15 cm}  \lambda(\mu_k^{\vee}) >0,  \\
  e^{s_{\mu_k}\lambda} -(1-e^{-\mu_k}) 
\big[ e^{\lambda+\mu_k}+ \cdots +
     e^{\lambda-\lambda(\mu_k^{\vee})\mu_k} \big]
 & { \rm if } \hspace{0,15 cm}  \lambda(\mu_k^{\vee}) <0  .
\end{array}\right.
$$

These relations hold since $s_{\mu_k}\lambda = \lambda - 
\lambda(\mu_k^{\vee}) \mu_k$.

\end{proof}

\begin{defi}

Let $\beta$ be a simple root, we define two $\mathbb{Z}$-linear maps
$T_{\beta}^0$ and $T_{\beta}^1$ from $R[H]$ to $R[H]$ by 
$$T_{\beta}^1(e^{\lambda})=e^{s_{\beta}\lambda},$$
$$
T_{\beta}^0(e^{\lambda})=\left\{ \begin{array}{ll}
 0  & { \rm if } \hspace{0,15 cm}  \lambda(\beta^{\vee})=0, \\ 
 e^{\lambda}+ e^{\lambda-\beta}+ \cdots +
     e^{\lambda-(\lambda(\beta^{\vee})-1)\beta}
 & { \rm if } \hspace{0,15 cm}  \lambda(\beta^{\vee}) >0,  \\ 
- e^{\lambda+\beta}- \cdots -  e^{\lambda-\lambda(\beta^{\vee})\beta} 
 & { \rm if } \hspace{0,15 cm}  \lambda(\beta^{\vee}) <0  ,
\end{array}\right.
$$
for all characters $e^{\lambda} \in X(H)$.

\end{defi}

\begin{rema}

The operators $T_{\beta}^0$ are called Demazure operators. Such
operators were first defined by Demazure in \cite{dema}.

\end{rema}

\begin{theo}  \label{linebundleBS}

For all integers $1 \leq k \leq N$ and 
all $\lambda \in \mathfrak{h}_{\mathbb{Z}}^*$, 
$$ [\mathcal{L}_{\lambda}^k]^H=\sum_{\epsilon \in \{0,1\}^k}
R_{\mu_1, \ldots ,\mu_k}^{\lambda,
  \epsilon}\mathcal{O}_{k,\epsilon}^H,$$
where $R_{\mu_1, \ldots, \mu_k}^{\lambda, \epsilon}=
T_{\mu_1}^{\epsilon_1}T_{\mu_2}^{\epsilon_2}   \cdots 
T_{\mu_k}^{\epsilon_k} (e^ {\lambda})$.

\end{theo}

\begin{proof}

We prove this theorem by induction on $k$.

For $k=1$, the theorem is a consequence of Theorem~\ref{induction} 
in the case $k=1$ and of the fact that $[\mathcal{L}_{\lambda}^0]^H=
e^{\lambda}$ for all $\lambda \in \mathfrak{h}_{\mathbb{Z}}^*$.

We assume that the relation is proved for $k-1$ ($1 \leq k-1 \leq 
N-1$) and for all weights $\lambda \in \mathfrak{h}_{\mathbb{Z}}^*$.  

We assume for example that $\lambda(\mu_k^{\vee})>0$. Then by 
Theorem \ref{induction}, we get
$$[\mathcal{L}_{\lambda}^k]^H=\sum_{\epsilon \in \{0,1\}^{k-1}} \! \!
\! \! \!
R_{\mu_1, \ldots ,\mu_{k-1}}^{s_{\mu_k}\lambda,\epsilon} 
\pi_k^*(\mathcal{O}_{k-1,\epsilon}^H) 
\, +\sum_{j=0}^{\lambda(\mu_k^{\vee})-1} \! \! \!
\sum_{\epsilon \in \{0,1\}^{k-1}} \! \! \! \mathcal{O}_{k,k}^H 
R_{\mu_1, \ldots ,\mu_{k-1}}^{\lambda-j\mu_k,\epsilon} 
\pi_k^*(\mathcal{O}_{k-1,\epsilon}^H). 
$$ 

Since $\pi_k$ is a smooth $H$-equivariant morphism between smooth 
$H$-varieties, for all $\epsilon \in \{0,1\}^{k-1}$, 
$$\pi_k^*(\mathcal{O}_{k-1,\epsilon}^H)=[\mathcal{O}_{\pi_k^{-1}
(\overline{\Gamma}_{\epsilon}^{k-1})}]^H=\mathcal{O}_{k,\tilde{\epsilon}^1}^H,$$ 
where $\tilde{\epsilon}^1=(\epsilon_1, \epsilon_2, \ldots,
\epsilon_{k-1}, 1)
\in \{0,1\}^k$. Moreover by Corollary \ref{product}, 
$$\mathcal{O}_{k,k}^H \pi_k^*(\mathcal{O}_{k-1,\epsilon}^H)=
\mathcal{O}_{k,k}^H \mathcal{O}_{k,\tilde{\epsilon}^1}^H
=\mathcal{O}_{k,\tilde{\epsilon}^0}^H,
$$
where $\tilde{\epsilon}^0=(\epsilon_1, \epsilon_2, \ldots, \epsilon_{k-1}, 0)
\in \{0,1\}^k$.

Then we get
$$
[\mathcal{L}_{\lambda}^k]^H=
\sum_{\tiny \begin{array}{c} \epsilon \in \{0,1\}^k \\ 
\epsilon_k=1 \end{array}}
R_{\mu_1, \ldots ,\mu_{k-1}}^{s_{\mu_k}\lambda,\epsilon}
\mathcal{O}_{k,\epsilon}^H
+ \sum_{ \tiny \begin{array}{c} \epsilon \in \{0,1\}^k \\ 
\epsilon_k=0 \end{array}} 
\sum_{j=0}^{\lambda(\mu_k^{\vee})-1} 
R_{\mu_1, \ldots ,\mu_{k-1}}^{\lambda-j\mu_k,\epsilon}
\mathcal{O}_{k,\epsilon}^H.
$$

We obtain the statement since by definition
$$R_{\mu_1, \ldots ,\mu_k}^{\lambda,\epsilon}=
\left\{ \begin{array}{ll}
R_{\mu_1, \ldots ,\mu_{k-1}}^{s_{\mu_k}\lambda,\epsilon}
 & { \rm if } \hspace{0,15 cm}  \epsilon_k=1, \\ 
 \sum_{j=0}^{\lambda(\mu_k^{\vee})-1} 
R_{\mu_1, \ldots ,\mu_{k-1}}^{\lambda-j\mu_k,\epsilon}
 & { \rm if } \hspace{0,15 cm}  \epsilon_k=0  .
\end{array}\right.
$$

The argument works in the same way in the cases
$\lambda(\mu_k^{\vee})=0$ and $\lambda(\mu_k^{\vee})<0$.

\end{proof}

\begin{exem}  \label{A2BS}

In type $A_2$ ($G=SL(3, \mathbb{C})$), we decompose
$[\mathcal{L}_{\rho_1}^{\Gamma}]^H$ in $K(H,\Gamma)$ where
$\Gamma=\Gamma(\alpha_2, \alpha_1, \alpha_2)$.

Since $\rho_1(\alpha_2^{\vee})=0$, we get
$$T_{\alpha_2}^0(e^{\rho_1})=0 \, \, { \rm and }
\, \, T_{\alpha_2}^1(e^{\rho_1})=e^{\rho_1} .$$ 

\medskip

Since $\rho_1(\alpha_1^{\vee})=1$, we find 
$$T_{\alpha_1}^0T_{\alpha_2}^1(e^{\rho_1})=e^{\rho_1} \, \, { \rm and
} \, \, T_{\alpha_1}^1T_{\alpha_2}^1(e^{\rho_1})=e^{\rho_1-\alpha_1}
 =e^{-\rho_1+\rho_2},$$
 since $\alpha_1 = 2 \rho_1 - \rho_2$.

\medskip
\medskip

Since $\rho_1(\alpha_2^{\vee})=0$ and  $(-\rho_1+\rho_2) 
(\alpha_2^{\vee})=1$, we find
$$T_{\alpha_2}^0T_{\alpha_1}^0T_{\alpha_2}^1(e^{\rho_1})=0, \,  
 \, T_{\alpha_2}^1 T_{\alpha_1}^0T_{\alpha_2}^1(e^{\rho_1})=
 e^{\rho_1}, $$
$$ T_{\alpha_2}^0T_{\alpha_1}^1T_{\alpha_2}^1(e^{\rho_1})=
 e^{-\rho_1 + \rho_2} , \, 
 T_{\alpha_2}^1 \, T_{\alpha_1}^1T_{\alpha_2}^1(e^{\rho_1})=
 e^{-\rho_1+ \rho_2-\alpha_2}=e^{-\rho_2},$$
since $\alpha_2 = -\rho_1 + 2\rho_2$.

\bigskip

Then Theorem \ref{linebundleBS} gives us the following relation in 
$K(H,\Gamma)$
$$[\mathcal{L}_{\rho_1}^{\Gamma}]^H=e^{-\rho_2} \mathcal{O}_{3,(1,1,1)}^H
+ e^{-\rho_1+\rho_2}   \mathcal{O}_{3,(0,1,1)}^H +
e^{\rho_1}\mathcal{O}_{3,(1,0,1)}^H. $$

\end{exem}

\begin{exem}  \label{G2BS}

In the case $G_2$, we decompose
$[\mathcal{L}_{\rho_2}^{\Gamma}]^H$ in $K(H,\Gamma)$ where we take 
$\Gamma=\Gamma(\alpha_1,\alpha_2, \alpha_1, \alpha_2)$.

Since $\rho_2=3\alpha_1+2 \alpha_2$, and 
$\rho_2(\alpha_2^{\vee})=1$, we get
$$T_{\alpha_2}^0(e^{\rho_2})=e^{3\alpha_1+2\alpha_2} \, \, { \rm and }
\, \, T_{\alpha_2}^1(e^{\rho_2})=e^{3\alpha_1+\alpha_2}.$$ 

\medskip

Since $(3\alpha_1+2\alpha_2)(\alpha_1^{\vee})=0$, 
and $(3\alpha_1+\alpha_2)(\alpha_1^{\vee})=3$, we find 
$$T_{\alpha_1}^0T_{\alpha_2}^0(e^{\rho_2})=0 \, ,\, \, 
T_{\alpha_1}^1T_{\alpha_2}^0(e^{\rho_2})=
e^{3\alpha_1+2\alpha_2}  \, , $$
$$ T_{\alpha_1}^0T_{\alpha_2}^1(e^{\rho_2})=e^{3 \alpha_1 +
 \alpha_2} + e^{2\alpha_1+\alpha_2} + e^{\alpha_1+\alpha_2}  \, \, 
{ \rm and} \, \, 
T_{\alpha_1}^1T_{\alpha_2}^1(e^{\rho_2})=e^{\alpha_2}.
$$

\medskip

Since $(3\alpha_1+2\alpha_2)(\alpha_2^{\vee})=1$,  
 $(3\alpha_1+\alpha_2)(\alpha_2^{\vee})=-1$, 
 $(2\alpha_1+\alpha_2)(\alpha_2^{\vee})=0$, 
 $(\alpha_1+\alpha_2)(\alpha_2^{\vee})=1$ and
 $\alpha_2(\alpha_2^{\vee})=2$, we find
$$T_{\alpha_2}^0T_{\alpha_1}^1T_{\alpha_2}^0(e^{\rho_2})=
e^{3\alpha_1+2\alpha_2}, \,  
 T_{\alpha_2}^1 T_{\alpha_1}^1T_{\alpha_2}^0(e^{\rho_2})=
 e^{3\alpha_1+\alpha_2}, $$
$$ T_{\alpha_2}^0T_{\alpha_1}^0T_{\alpha_2}^1(e^{ \rho_2})=
 -e^{3\alpha_1+2\alpha_2}+e^{\alpha_1+\alpha_2} , \, 
 T_{\alpha_2}^1 \, T_{\alpha_1}^0T_{\alpha_2}^1(e^{ \rho_2})=
 e^{3\alpha_1+2\alpha_2}+e^{2\alpha_1+\alpha_2}
 +e^{\alpha_1},$$
$$ T_{\alpha_2}^0T_{\alpha_1}^1T_{\alpha_2}^1(e^{ \rho_2})=
 e^{\alpha_2}+1 , \, 
 T_{\alpha_2}^1 \, T_{\alpha_1}^1T_{\alpha_2}^1(e^{ \rho_2})=
 e^{-\alpha_2}.$$

\medskip

Since $(3\alpha_1+2\alpha_2)(\alpha_1^{\vee})=0$,  
 $(3\alpha_1+\alpha_2)(\alpha_1^{\vee})=3$, 
 $(\alpha_1+\alpha_2)(\alpha_1^{\vee})=-1$, 
 $(2\alpha_1+\alpha_2)(\alpha_1^{\vee})=1$,
 $\alpha_1(\alpha_1^{\vee})=2$, and
 $\alpha_2(\alpha_1^{\vee})=-3$, we find

$$T_{\alpha_1}^0T_{\alpha_2}^0T_{\alpha_1}^1
T_{\alpha_2}^0(e^{\rho_2})=0, \,  
 T_{\alpha_1}^1T_{\alpha_2}^0T_{\alpha_1}^1T_{\alpha_2}^0
 (e^{\rho_2})=e^{3\alpha_1+2\alpha_2}, 
 $$

 $$T_{\alpha_1}^0T_{\alpha_2}^1 T_{\alpha_1}^1
 T_{\alpha_2}^0(e^{\rho_2})=e^{3\alpha_1+\alpha_2}
 +e^{2\alpha_1+\alpha_2}+e^{\alpha_1+\alpha_2}, \, 
 T_{\alpha_1}^1T_{\alpha_2}^1 T_{\alpha_1}^1T_{\alpha_2}^0
 (e^{\rho_2})=e^{\alpha_2}, 
 $$

$$T_{\alpha_1}^0 T_{\alpha_2}^0T_{\alpha_1}^0
T_{\alpha_2}^1(e^{ \rho_2})= -e^{2\alpha_1+\alpha_2}, \,
T_{\alpha_1}^1T_{\alpha_2}^0T_{\alpha_1}^0T_{\alpha_2}^1
(e^{ \rho_2})= -e^{3\alpha_1+2\alpha_2}+e^{2\alpha_1+\alpha_2} ,
 $$ 

$$ T_{\alpha_1}^0T_{\alpha_2}^1 \, T_{\alpha_1}^0
T_{\alpha_2}^1(e^{ \rho_2})= e^{2\alpha_1+\alpha_2}+
e^{\alpha_1} + 1,  \, 
 T_{\alpha_1}^1T_{\alpha_2}^1 \, T_{\alpha_1}^0T_{\alpha_2}^1
 (e^{ \rho_2})=e^{3\alpha_1+2\alpha_2}+e^{\alpha_1+\alpha_2}
 +e^{-\alpha_1}, 
 $$

$$ T_{\alpha_1}^0T_{\alpha_2}^0T_{\alpha_1}^1T_{\alpha_2}^1
(e^{ \rho_2})= -e^{\alpha_1+\alpha_2}-e^{2\alpha_1+\alpha_2}
 -e^{3\alpha_1 +\alpha_2} , \, 
 T_{\alpha_1}^1T_{\alpha_2}^0T_{\alpha_1}^1T_{\alpha_2}^1
(e^{ \rho_2})= e^{3\alpha_1+\alpha_2}+1 ,
$$

 $$
 T_{\alpha_1}^0T_{\alpha_2}^1 \, T_{\alpha_1}^1T_{\alpha_2}^1
 (e^{ \rho_2})= e^{-\alpha_2}+e^{-\alpha_1-\alpha_2}
 +e^{-2\alpha_1 -\alpha_2} , \, 
  T_{\alpha_1}^1T_{\alpha_2}^1 \, T_{\alpha_1}^1T_{\alpha_2}^1
  (e^{ \rho_2})= e^{-3\alpha_1-\alpha_2}.
$$

\bigskip

Then Theorem \ref{linebundleBS} gives us the following relation in 
$K(H,\Gamma)$
$$[\mathcal{L}_{\rho_2}^{\Gamma}]^H=
 e^{-3\alpha_1-\alpha_2} \mathcal{O}_{4,(1,1,1,1)}^H
+ (e^{-\alpha_2}+e^{-\alpha_1-\alpha_2}
 +e^{-2\alpha_1 -\alpha_2})\mathcal{O}_{4,(0,1,1,1)}^H $$
$$+
(e^{3\alpha_1+\alpha_2}+1)\mathcal{O}_{4,(1,0,1,1)}^H  +
 ( -e^{\alpha_1+\alpha_2}-e^{2\alpha_1+\alpha_2}
 -e^{3\alpha_1 +\alpha_2}) \mathcal{O}_{4,(0,0,1,1)}^H
 $$
$$+   (e^{3\alpha_1+2\alpha_2}+e^{\alpha_1+\alpha_2}
 +e^{-\alpha_1}) \mathcal{O}_{4,(1,1,0,1)}^H +
(e^{2\alpha_1+\alpha_2}+
e^{\alpha_1} + 1)\mathcal{O}_{4,(0,1,0,1)}^H $$
$$  +(-e^{3\alpha_1+2\alpha_2}+e^{2\alpha_1+\alpha_2})
\mathcal{O}_{4,(1,0,0,1)}^H
 -e^{2\alpha_1+\alpha_2}\mathcal{O}_{4,(0,0,0,1)}^H $$
 $$+   e^{\alpha_2}\mathcal{O}_{4,(1,1,1,0)}^H  +
(e^{3\alpha_1+\alpha_2}
 +e^{2\alpha_1+\alpha_2}+e^{\alpha_1+\alpha_2})
  \mathcal{O}_{4,(0,1,1,0)}^H
+e^{3\alpha_1+2\alpha_2}  \mathcal{O}_{4,(1,0,1,0)}^H  . $$

\end{exem}

\section{Equivariant $K$-theory of flag varieties} \label{sectionvd}

\subsection{Definitions}

Since $X$ is a complex irreducible smooth $H$-variety, we denote by
$K(H,X)$ its $H$-equivariant $K$-theory.

For $w \in W$, we set  $\mathcal{O}_w^H = [\mathcal{O}_{\overline{X}_w}]^H \in
K(H,X)$.

Since $X = \coprod_{w \in W} X_w$ is a cell decomposition of $X$, the
family  $\{\mathcal{O}_w^H \}_{w \in W}$ is a
$R[H]$-basis of the module $K(H,X)$. Moreover, $X^H$
is finite and isomorphic to $W$. Thus we have the following
proposition.

\begin{prop} \label{propbasektheorie}

\indent

\begin{enumerate}

\item[$(i)$] The $H$-equivariant $K$-theory of $X^H$ can be
  identified with $F(W;R[H])$.

\item[$(ii)$]  $ K(H,X)=\oplus_{w \in W} R[H] \mathcal{O}_w^H$.  

\item[$(iii)$] The restriction to fixed points $i_{H}^*$ :
  $K(H,X) \rightarrow F(W;R[H])$ is injective.

\end{enumerate}

\end{prop}

For all $\lambda \in \mathfrak{h}_{\mathbb{Z}}^*$, 
we denote by $\mathcal{L}_{\lambda}^{X}$ the $B$-equivariant line
bundle over $X$ defined as the space of orbits of $B$ acting on   
$G \times \mathbb{C}$ by
$$(g, v) b = (gb, e^{\lambda}(b^{-1})v), b \in B, \, g \in G, \, v \in 
\mathbb{C}.$$

\subsection{Link with Bott-Samelson varieties}

Let $\mu_{1}, \ldots , \mu_{N}$ be a sequence of $N$ simple roots. 
We set $\Gamma=\Gamma(\mu_{1}, \ldots, \mu_{N})$ and we define
an $H$-equivariant map $g$ from $\Gamma$ to $X$ by multiplication  
$$g([g_{1},\ldots ,g_{N}]) = g_{1}\times \cdots   \times g_{N} \;
\,[B].$$

\medskip

For all $\epsilon \in \{ 0, 1\}^N$, we define
$$\displaystyle{\underline{v}(\epsilon) = 
\prod_{ {\tiny \begin{array}{cc}  1\leq j \leq N, \\ 
  \epsilon_j =1 \end{array}} } 
 \!\!\!\!\!\!  \underline{s_{\mu_{j}}} \in \underline{W}}.$$ 
For all $\epsilon \in \{ 0, 1\}^N$, we set $\mathcal{O}_{\epsilon}^H
= \mathcal{O}_{N,\epsilon}^H \in K(H,\Gamma)$ (see Section
~\ref{basisBS}). We have the following theorem.

\begin{theo}  \label{g*}

For all $\epsilon \in \{ 0, 1\}^N$,
$$g_*(\mathcal{O}_{\epsilon}^H)=
\mathcal{O}_{T^{-1}(\underline{v}(\epsilon))}^H,
$$
where $T$ is defined in Section \ref{monoid}.

\end{theo}

\begin{proof}

By Theorem~8.1.13 and Corollary~8.2.3 of \cite{livrekumar}, 
the image of $g$ is a Schubert variety $\overline{X}_w$, where $w \in W$,
and $g_*(\mathcal{O}_{\epsilon}^H)=\mathcal{O}_w^H$.
 By Lemma~2.3 of \cite{mw3},  $w=T^{-1}(\underline{v}(\epsilon))$.

\end{proof}

\begin{exem}  \label{A2g*}

In the case $A_2$, if we take $\Gamma=~\Gamma(\alpha_2, \alpha_1,
\alpha_2)$, we get the following relations : 
$$
\begin{array}{lcl}
\vspace{0.1 cm}
g_*(\mathcal{O}_{(0,0,0)}^H) & = & \mathcal{O}_1^H  \\
\vspace{0.1 cm}
g_*(\mathcal{O}_{(0,0,1)}^H) & = & \mathcal{O}_{s_2}^H  \\
\vspace{0.1 cm}
g_*(\mathcal{O}_{(0,1,0)}^H) & = & \mathcal{O}_{s_1}^H  \\
\vspace{0.1 cm}
g_*(\mathcal{O}_{(1,0,0)}^H) & = & \mathcal{O}_{s_2}^H  \\
\vspace{0.1 cm}
g_*(\mathcal{O}_{(1,0,1)}^H) & = & \mathcal{O}_{s_2}^H \\
\vspace{0.1 cm}
g_*(\mathcal{O}_{(0,1,1)}^H) & = & \mathcal{O}_{s_1s_2}^H \\
\vspace{0.1 cm}
g_*(\mathcal{O}_{(1,1,0)}^H) & = & \mathcal{O}_{s_2s_1}^H  \\
\vspace{0.1 cm}
g_*(\mathcal{O}_{(1,1,1)}^H) & = & \mathcal{O}_{s_2s_1s_2}^H.

\end{array}
$$

\end{exem}

\begin{lemm} \label{g*linebundles}
Let $w=s_{\mu_1}  \cdots s_{\mu_N}$ be a reduced decomposition of $w
\in W$ (N=l(w)). For all $\lambda \in \mathfrak{h}_{\mathbb{Z}}^*$,
$$g_*([\mathcal{L}_{ \lambda }^{\Gamma}]^H)=
[\mathcal{L}_{ \lambda }^X]^H \times \mathcal{O}_w^H.
$$

\end{lemm}

\begin{proof}
Since $ \mathcal{O}_{(\bold{1})}^H=1 \in K(H, \Gamma)$ where
$(\bold{1})=(1,1,\ldots,1) \in \{0,1\}^N$,
and for all $\lambda \in \mathfrak{h}_{\mathbb{Z}}^*$, $g^*([\mathcal{L}_{\lambda}^X]^H)=
[\mathcal{L}_{\lambda}^{\Gamma}]^H$, 
we have 
$$g_*([\mathcal{L}_{ \lambda }^{\Gamma}]^H)
=g_*(g^*([\mathcal{L}_{\lambda}^X]^H) \times \mathcal{O}_{(\bold{1})}^H)
=[\mathcal{L}_{\lambda}^X]^H \times g_*(\mathcal{O}_{(\bold{1})}^H)
= [\mathcal{L}_{\lambda}^X]^H \times \mathcal{O}_w^H,$$
where the last equality is a consequence of Theorem~\ref{g*}.

\end{proof}

\subsection{A Chevalley formula}

Theorems~\ref{linebundleBS}, \ref{g*} and Lemma 
\ref{g*linebundles} give us the following theorem.

\begin{theo}  \label{main}
Let $w=s_{\mu_1}  \cdots s_{\mu_N}  $ be a reduced 
decomposition of $w \in W$. For all $\lambda \in 
\mathfrak{h}_{\mathbb{Z}}^*$,
$$
[\mathcal{L}_{\lambda}^X]^H \times \mathcal{O}_w^H =
\sum_{\epsilon \in \{0,1\}^N} R_{\mu_1, \ldots, \mu_N}^{\lambda,
  \epsilon} \mathcal{O}_{T^{-1}(\underline{v}(\epsilon))}^H.
$$

\end{theo}

\begin{exem} \label{A2H}

In the case $A_2$, if we take $w=s_2s_1s_2$ and
$\lambda = \rho_1$, Examples~\ref{A2BS}, \ref{A2g*} and 
Theorem~\ref{main} give us the following relation in $K(H,X)$
$$ [\mathcal{L}_{\rho_1}^X]^H \times \mathcal{O}_{s_2s_1s_2}^H= 
\mathcal{L}_{\rho_1}^X= e^{-\rho_2} \mathcal{O}_{s_2s_1s_2}^H
+   e^{-\rho_1+\rho_2} \mathcal{O}_{s_1s_2}^H + e^{\rho_1}
\mathcal{O}_{s_2}^H . $$ 
 
\end{exem}

\begin{exem} \label{G2H}

In the case $G_2$, if we take $w=s_1s_2s_1s_2$ and
$\lambda = \rho_2$, Example~\ref{G2BS} and 
Theorem~\ref{main} give us the following relation in $K(H,X)$
$$ [\mathcal{L}_{\rho_2}^X]^H \times \mathcal{O}_{s_1s_2s_1s_2}^H 
= e^{-3\alpha_1-\alpha_2} \mathcal{O}_{s_1s_2s_1s_2}^H
+ (e^{-\alpha_2}+e^{-\alpha_1-\alpha_2}
 +e^{-2\alpha_1 -\alpha_2})\mathcal{O}_{s_2s_1s_2}^H $$
$$+
(e^{3\alpha_1+\alpha_2}+1)\mathcal{O}_{s_1s_2}^H  +
 ( -e^{\alpha_1+\alpha_2}-e^{2\alpha_1+\alpha_2}
 -e^{3\alpha_1 +\alpha_2}) \mathcal{O}_{s_1s_2}^H
 $$
$$+   (e^{3\alpha_1+2\alpha_2}+e^{\alpha_1+\alpha_2}
 +e^{-\alpha_1}) \mathcal{O}_{s_1s_2}^H +
(e^{2\alpha_1+\alpha_2}+
e^{\alpha_1} + 1)\mathcal{O}_{s_2}^H $$
$$  +(-e^{3\alpha_1+2\alpha_2}+e^{2\alpha_1+\alpha_2})
\mathcal{O}_{s_1s_2}^H
 -e^{2\alpha_1+\alpha_2}\mathcal{O}_{s_2}^H $$
 $$+   e^{\alpha_2}\mathcal{O}_{s_1s_2s_1}^H  +
(e^{3\alpha_1+\alpha_2}
 +e^{2\alpha_1+\alpha_2}+e^{\alpha_1+\alpha_2})
  \mathcal{O}_{s_2s_1}^H
+e^{3\alpha_1+2\alpha_2}  \mathcal{O}_{s_1}^H   $$
$$ =e^{-3\alpha_1-\alpha_2} \mathcal{O}_{s_1s_2s_1s_2}^H
+  (e^{-\alpha_2}+e^{-\alpha_1-\alpha_2}
 +e^{-2\alpha_1 -\alpha_2})\mathcal{O}_{s_2s_1s_2}^H $$
$$ +   e^{\alpha_2}\mathcal{O}_{s_1s_2s_1}^H
+ (e^{3\alpha_1+\alpha_2}
+e^{2\alpha_1+\alpha_2}+e^{\alpha_1+\alpha_2}) 
\mathcal{O}_{s_2s_1}^H+(e^{-\alpha_1}+1) \mathcal{O}_{s_1s_2}^H 
 $$
 $$ +(e^{\alpha_1} + 1)\mathcal{O}_{s_2}^H
+e^{3\alpha_1+2\alpha_2}  \mathcal{O}_{s_1}^H.   $$

\end{exem}

\bigskip

\begin{rema}

Since $\rho_2$ is a dominant weight (i.e. $\rho_2(\alpha_i^{\vee})
\geq 0$ for all simple roots $\alpha_i$), we know that we must find 
positive coefficients (i.e. a linear combination of characters with 
positive coefficients, see \cite{mathieu}). Unfortunately, our formula
is not positive. In this example, we see that we can find negative
terms and then cancellations can occur. The formulas given by Pittie
and Ram in \cite{pittie} and Littelmann and Seshadri in \cite{littel}
are positive. 

\end{rema}

\section{Ordinary $K$-theory}  \label{sectionordinary}

We denote by $\psi$ the forgetful map $K(H,X) \rightarrow 
K(X)$, where $K(X) \simeq K^0(X) \simeq K_0(X)$ is the ordinary
$K$-theory of $X$, and by $ev :
R[H] \rightarrow \mathbb{Z}$ the $\mathbb{Z}$-linear map 
defined by $$ev(e^{\alpha})=1 {\rm \,\,\, for \, all \, character
  \,\,} e^{\alpha} \in  
X(H).$$ 

For all $w \in W$, we denote by $\mathcal{O}_w \in 
K(X)$ the class of $\mathcal{O}_{\overline{X}_w}$ in 
$K(X)$, and for all $\lambda \in \mathfrak{h}_{\mathbb{Z}}^*$, we 
set $[\mathcal{L}_{\lambda}^X]=\psi([\mathcal{L}_{\lambda}^X]^H) \in 
K(X)$. Since $\psi(\mathcal{O}_w^H)=\mathcal{O}_w$, 
$\psi(e^{\alpha})=ev(e^{\alpha})$, and $\psi$ is a ring homomorphism, 
Theorem \ref{main} gives us the following theorem.

\begin{theo} \label{mainordinary}

Let $w=s_{\mu_1}  \cdots s_{\mu_N}  $ be a reduced 
decomposition of $w \in W$. For all $\lambda \in 
\mathfrak{h}_{\mathbb{Z}}^*$,
$$
[\mathcal{L}_{\lambda}^X] \times \mathcal{O}_w =
\sum_{\epsilon \in \{0,1\}^N} ev(R_{\mu_1, \ldots, \mu_N}^{\lambda,
  \epsilon}) \mathcal{O}_{T^{-1}(\underline{v}(\epsilon))}.
$$

\end{theo}

\begin{exem} \label{A2ordinary}

In the case $A_2$, if we take
$w=s_2s_1s_2$ and $\lambda = \rho_1$, Example~\ref{A2H} and 
Theorem~\ref{mainordinary} give us the following relation in $K(X)$
$$[\mathcal{L}_{\rho_1}^X] \times \mathcal{O}_{s_2s_1s_2}= 
[\mathcal{L}_{\rho_1}^X]= \mathcal{O}_{s_2s_1s_2}
+\mathcal{O}_{s_1s_2} + \mathcal{O}_{s_2} . $$

\end{exem}

\begin{exem} \label{G2ordinary}

In the case $G_2$, if we take
$w=s_1s_2s_1s_2$ and $\lambda = \rho_2$, Example~\ref{G2H} and 
Theorem~\ref{mainordinary} give us the following relation in $K(X)$
$$[\mathcal{L}_{\rho_2}^X] \times \mathcal{O}_{s_1s_2s_1s_2}= 
 \mathcal{O}_{s_1s_2s_1s_2} + 3\mathcal{O}_{s_2s_1s_2} 
 +\mathcal{O}_{s_1s_2s_1} +  3 \mathcal{O}_{s_2s_1} +
 2 \mathcal{O}_{s_1s_2} +2 \mathcal{O}_{s_2}
+ \mathcal{O}_{s_1}  . $$
 
This example was computed by Pittie and Ram in \cite{pittie} by using
LS-paths.

\end{exem}

   \bibliography{chevalley}
   \bibliographystyle{amsplain}

  \end{document}